\documentclass[a4paper,12pt,leqno]{article}
\usepackage{amsmath,amsfonts,latexsym}
\usepackage{amssymb}
\usepackage{graphics}

\newcommand{\nc}{\newcommand}
\nc{\tr}{{\triangle}} \nc{\vth}{{\vartheta}} \nc{\bt}{{\beta}}
\nc{\dl}{{\delta}} \nc{\Dl}{{\Delta}} \nc{\p}{{\psi}}
\nc{\gm}{{\gamma}} \nc{\Gm}{{\Gamma}} \nc{\sg}{{\sigma}}
\nc{\ve}{{\varepsilon}} \nc{\ch}{{\cal H}} \nc{\cf}{{\cal F}}
\nc{\cp}{{\cal P}} \nc{\td}{\tilde}

\newtheorem{rem}{Remark}[section]

\newtheorem{cor}{Corollary}[section]
\newtheorem{Lemma}{Lemma}[section]

\begin{document}

%***************************************************************

%***************************************************************************

\title
{Rate of Convergence in   Recursive Parameter Estimation procedures.}
   \author{Teo Sharia}
\date{}
\maketitle
\begin{center}
{\it
Department of Mathematics \\Royal Holloway,  University of London\\
Egham, Surrey TW20 0EX \\ e-mail: t.sharia@rhul.ac.uk }
\end{center}
\begin{abstract}
 We consider
 estimation procedures which are recursive in the sense that each successive
 estimator is obtained from the previous one by a simple adjustment. We study
         rate of convergence of recursive estimation
       procedures for the general statistical model.
\end{abstract}
\begin{center}
\end{center}

                       %S E C T I O N 1
\section{Introduction}
Let $X_1, \dots, X_n$ be  random variables, with a
joint distribution depending on a real  unknown parameter $\theta$. Then an
$M$-estimator of $\theta$ is defined as  a solution of the
estimating equation
%%%%%%%%%%%%%%%%%%%%%%%%%%%%%%%%%%%%%%%%%%%%%%%%%%%%
%                  \eqref{esteqg}
%%%%%%%%%%%%%%%%%%%%%%%%%%%%%%%%%%%%%%%%%%%%%%%%%%%%%%
\begin{equation}\label{esteqg}
\sum_{i=1}^n \psi_i(v)=0,
\end{equation}
where   $\psi_i(v)=\psi_i(X_{i-k}^i; v)$  ~($i=1,2,\dots, n$)
are suitably chosen functions and
$X_{i-k}^i=
(X_{i-k}, \dots,X_i )$ is the a vector of past and present observations
 at step (time) $i$.  For instance, if $X_i$'s are
observations from a discrete time Markov process, then one can
assume that   $k=1$.
If observations are i.i.d., then we take $k=0$ so that
$\psi_i(v)=\psi_i(X_{i}; v).$
In general, if no restrictions are made on
the dependence structure of the process $X_i$, one may need to
consider $\psi$-functions depending on the vector of all past and
present observations of the process (that is, $k=i-1$). If the
conditional probability density function (or probability function)
of the observation $X_i,$ given $X_{i-k}, \dots,X_{i-1},$ is
$f_i(x,\theta)=f_i(x,\theta|X_{i-k}, \dots,X_{i-1})$,  then one
can obtain a   MLE (maximum likelihood estimator)  on choosing
$\psi_i(v)=f'_i(X_i,v)/f_i(X_i,v).$ Besides MLEs, the class of
$M$-estimators includes estimators
 with special properties such as {robustness}.
Under certain regularity and ergodicity conditions it can be
proved that there exists a consistent sequence of solutions of
\eqref{esteqg} which has  the property of local asymptotic
linearity (See e.g., Serfling (1980),
Huber (1981), Lehman (1983).  A comprehensive bibliography can be
found in Launer and Wilkinson (1979), Hampel at al (1986), Rieder
(1994), and Jure$\check{c}$kov$\acute{a}$ and Sen (1996).)

 If   $\psi$-functions are
 nonlinear, it is rather difficult to
 work with  the corresponding estimating equations.
  In this paper we consider
 estimation procedures which are recursive in the sense that each successive
 estimator is obtained from the previous one by a simple adjustment.
   In particular, we consider  a class of
estimators
%%%%%%%%%%%%%%%%%%%%%%%%%%%%%%%%%%%%%%%%%%%%%%%%%%
%                   \eqref{rec1}
%%%%%%%%%%%%%%%%%%%%%%%%%%%%%%%%%%%%%%%%%%%%%%%%%
\begin{equation} \label{rec1}
\hat\theta_n=\hat\theta_{n-1}+
{\Gamma_n^{-1}(\hat\theta_{n-1})}\psi_n(\hat\theta_{n-1}),
~~~~~~~~~n\ge 1,
\end{equation}
where $\psi_n$ is a suitably chosen vector process, $\Gamma_n$  is
a (possibly random) normalizing matrix process and
$\hat\theta_0\in {\mathbb{R}}^m$ is some initial point.
(See the introduction in  Sharia (2006) for a detailed discussion and a heuristic justification
of this estimation procedure.)

 In  i.i.d. models,
 estimating procedures similar to \eqref{rec1} have been studied by a number of
authors using methods of stochastic approximation theory (see,
e.g., Khas'minskii  and Nevelson (1972), Fabian (1978),  Ljung and Soderstrom
(1987), Ljung, Pflug and Walk (1992),  and references therein).
Some work has been done for  non i.i.d. models as well. In
particular,  Englund,   Holst, and  Ruppert (1989) give an
asymptotic representation results for  certain type of  $X_n$
processes. In Sharia (1998) theoretical results on convergence,
rate of convergence and the asymptotic representation are given
under certain regularity and ergodicity assumptions on the model,
in the one-dimensional  case with
$\psi_n(x,\theta)=\frac{\partial}{\partial\theta} \mbox{log}
f_n(x,\theta)$ (see also Campbell (1982),
%Sharia (1992),
 Sharia (1997),
  Lazrieva and  Toronjadze  (1987)).

In Sharia (2006), imposing ``global'' restrictions on the processes
$\psi$ and $\Gamma$, we study   ``global'' convergence of the
recursive estimators \eqref{rec1}, that is,  convergence  for an arbitrary
starting point $\hat\theta_0$.
In the present paper,  we present results on  rate of the
convergence and demonstrate the use of these
results on some examples.

                       %S E C T I O N 2

\section {Notation and preliminaries}
\setcounter{equation}{0}

Let   $X_t,\;\; t=1,2,\dots ,$
 be  observations taking values in a  measurable space
$({\bf X},{\cal B}({\bf X}))$  equipped with  a $\sigma$-finite
measure $\mu.$ Suppose  that the distribution of the process $X_t$
depends on an unknown parameter $\theta \in \Theta,$ where
$\Theta$ is an open subset of the $m$-dimensional Euclidean space
$\mathbb{R}^m$. Suppose also  that for each  $t=1,2,\dots$, there
exists a  regular conditional probability density of $X_t$ given
values of past observations of  $X_{t-1},\dots , X_2, X_1$, which will be
denoted by
$$
f_ t(\theta, x_t \mid x_1^{t-1})= f_ t(\theta, x_t \mid
x_{t-1},\dots ,x_1),
$$
where $ f_ 1(\theta, x_1 \mid x_1^0)= f_ 1(\theta, x_1) $ is the
probability  density  of the random variable $X_1.$ Without loss of
generality we assume that
all random variables are defined on a probability space $
 (\Omega, \cf)
 $
and denote by $ \left\{P^\theta, \; \theta\in \Theta\right\} $ the
family of the corresponding distributions on $
 (\Omega, \cf).
 $

  Let  $\cf_t=\sigma (X_1,\dots ,X_t)$ be the $\sigma$-field
generated by the random variables $ X_1,\dots ,X_t.$ By
$\left(\mathbb{R}^m, {\cal B}( \mathbb{R}^m ) \right)$ we denote
the $m$-dimensional Euclidean space with the Borel
$\sigma$-algebra ${\cal B}( \mathbb{R}^m )$. Transposition of
matrices and vectors is denoted by $T$. By $(u,v)$ we denote the
standard scalar product of $u,v \in \mathbb{R}^m,$ that is,
$(u,v)=u^Tv.$

Suppose that $h$ is a real valued function  defined on
 $\Theta\subset  {{\mathbb{R}}}^m$.  We denote by $\dot h(\theta)$  the row-vector
 of partial derivatives
of $h(\theta)$ with respect to the components of $\theta$, that
is,
 $$
 \dot h(\theta)=\left(\frac{{\partial}}{{\partial} \theta^1} h(\theta), \dots,
 \frac{{\partial}}{{\partial} \theta^m} h(\theta)\right).
 $$

 If for each  $t=1,2,\dots$, the derivative
$ \dot f_t(\theta, x_t \mid x_1^{t-1})$ w.r.t. $\theta$
  exists, then we
can define the function
$$
l_t(\theta, x_t \mid x_1^{t-1})=\frac 1 {f_ t(\theta, x_t \mid
x_1^{t-1})} \dot f_t^T(\theta, x_t \mid x_1^{t-1})
$$
with the convention $0/0=0$.

The {\it one step   conditional  Fisher information matrix} for
$t=1,2,\dots$ is defined as
$$
i_t(\theta\mid x_1^{t-1}) =\int l_t(\theta,z\mid x_1^{t-1})
l^{T}_t(\theta,z\mid x_1^{t-1})f_t(\theta,z\mid x_1^{t-1})\mu
(dz).
$$
We shall use the  notation
$$
f_t(\theta)=f_t(\theta, X_t \mid X_1^{t-1}), \;\;\; \;\;\;
l_t(\theta)=l_t(\theta, X_t\mid X_1^{t-1}),
$$
$$
i_t(\theta)=i_t(\theta\mid X_1^{t-1}).
$$
Note that  the process $i_t(\theta)$  is {\it ``predictable''},
that is, the random variable $i_t(\theta),$ is $\cf_{t-1}$
measurable for each $t\ge 1.$

 Note also that by definition,
$i_t(\theta)$ is   a version of the conditional expectation w.r.t.
${\cal{F}}_{t-1},$  that is,
$$
i_t(\theta)= E_\theta\left\{l_t(\theta) l^{T}_t(\theta) \mid
{\cal{F}}_{t-1}\right\}.
$$
Everywhere in the present work conditional expectations  are meant
to be
 calculated as integrals w.r.t. the conditional probability densities.

The {\it conditional  Fisher information} at time $t$ is
$$
I_t(\theta)=\sum_{s=1}^t i_s(\theta),  \;\;\;\;\;\;\;\;\;
t=1,2,\dots.
$$
If the $X_t$'s are independent random variables,  $I_t(\theta)$
reduces to the standard Fisher information matrix.  Sometimes
$I_t(\theta)$ is referred as the incremental expected Fisher
information. Detailed discussion of this concept and related work
appears in
 Barndorff-Nielsen and  Sorensen (1994), and Prakasa-Rao  (1999)  Ch.3.

We say that  ${\bf \psi}= \{\psi_t(\theta, x_t, x_{t-1}, \dots,
x_1)\}_{t\ge 1}$ is a sequence of estimating functions and write
$\bf \psi \in \Psi$, if for each ${t\ge 1},$ $ \psi_t(\theta, x_t,
x_{t-1}, \dots, x_1)  :
 \Theta \times {\bf X}^t \;\;\to \;\;{\mathbb{R}}^m
$ is a   Borel function.

Note that  $
\{l_t(\theta, x_t \mid x_1^{t-1})\}_{t\ge 1}\in {\bf\Psi}$ and
a ML recursive procedure  is given by
$$
\hat \theta_t=\hat \theta_{t-1}+I_t^{-1}(\hat
\theta_{t-1})l_t(\hat \theta_{t-1}), \qquad  t\ge1.
$$

\vskip+0.5cm

 {\bf Convention} {\it  Everywhere in the present work
 $\theta\in \mathbb{R}^m $ is an arbitrary but fixed value
of the parameter.  Convergence and all relations between random
variables are meant with probability one w.r.t. the measure
$P^\theta$ unless specified otherwise. A sequence of random
variables $(\xi_t)_{t\ge1}$ has some property eventually if for
every $\omega$ in a set $\Omega^\theta$ of $P^\theta$ probability
1,
 $\xi_t$  has this property for all $t$ greater than some
$t_0(\omega)<\infty$.} \vskip+0.5cm

%%%%%%%%%%%%%%%%%%%%%%%%%%%%%%%%%%%%%%%%%%%%%%%%%%%%%%%%%%%%%%%%%%%%%%%%%%
%
%                    R A T E   O F   C O N V E R G E N C E
%
%%%%%%%%%%%%%%%%%%%%%%%%%%%%%%%%%%%%%%%%%%%%%%%%%%%%%%%%%%%%%%%%%%%%%%%%%%

\section{Main results}
\setcounter{equation}{0}

Suppose   that $\bf \psi \in \Psi$ and    $\Gm_t(\theta)$,  for
each $\theta\in \mathbb{R}^m$, is  a predictable $m\times m$ matrix process
with $ \mbox{det} ~\Gm_t(\theta)\neq0$, $t\ge 1$.
 Consider  the estimator $\hat \theta_t$  defined by
                       %(rec2)
\begin{equation}\label{rec2}
\hat \theta_t=\hat \theta_{t-1}+\Gm_t^{-1}(\hat
\theta_{t-1})\p_t(\hat \theta_{t-1}), \qquad  t\ge1,
\end{equation}
where  $\hat \theta_0\in\mathbb{R}^m $ is arbitrary initial point.

Let $\theta\in \mathbb{R}^m $ be an arbitrary but fixed value of
the parameter  and  for any  $u\in \mathbb{R}^m$ define
 $$
b_t(\theta,u)=E_\theta \left\{\p_t(\theta+u)\mid
{\cf}_{t-1}\right\}.
$$
%=E_\theta \left\{\p_t(\theta+u)-\p_t(\theta)\mid{\cf}_{t-1}\right\}.

%%%%%%%%%%%%%%%%%%%%%%%%%%%%%%%%%%%%%%%%%%%%%%%%%%%%%%%%%%%%%%%%%%%%%%%%%
%
%                          L e m m a   3.1
%
%%%%%%%%%%%%%%%%%%%%%%%%%%%%%%%%%%%%%%%%%%%%%%%%%%%%%%%%%%%%%%%%%%%%%%%%%%

\begin{Lemma} Let $\{C_t(\theta)\}$ be a symmetric predictable  $m\times m$ matrix process
such that $C_t(\theta)$ is non-negative  definite  for $t=1,2,\dots$.
Denote $\Dl_t=\hat\theta_t-\theta,$ ~ $V_t(u)=(C_t(\theta) u,u)$ and $ \tr
V_t(u)=V_t(u)-V_{t-1}(u).$ Suppose that
%                       (Kt)
\begin{eqnarray}\label{Kt}
\sum_{t=1}^{\infty}\left(1+ V_{t-1}(\Dl_{t-1})\right)^{-1}\left[
{\cal K}_t(\theta)\right]^+< \infty, \qquad
 P^\theta\mbox{-a.s.},
 \end{eqnarray}
where
%                      (Kt2)
\begin{eqnarray}\label{Kt2}
{\cal K}_t(\theta) = \tr V_t(\Dl_{t-1}) +2\left(C_t(\theta)
\Dl_{t-1}, \Gm_t^{-1}(\theta+\Dl_{t-1})
b_t(\theta,\Dl_{t-1})\right)\\
  + E_\theta \left\{\left[\Gm_t^{-1}(\theta+\Dl_{t-1})
\p_t(\theta+\Dl_{t-1})\right]^T C_t(\theta)
\Gm_t^{-1}(\theta+\Dl_{t-1}) \p_t(\theta+\Dl_{t-1})\mid
{\cf}_{t-1}\right\}. \nonumber
 \end{eqnarray}
Then $
 V_t(\Dl_{t})
$
 converges
 ($P^\theta$-a.s.)  to a finite limit.
\end{Lemma}
%%%%%%%%%%%%%%%%%%%%%%%%%%%%%%%%%%%%%%%%%%%%%%%%%%%%%%%%%%%%%%%%%%%%%%%%%%
%                  P R O O F   OF       L E M M A   3.1
%%%%%%%%%%%%%%%%%%%%%%%%%%%%%%%%%%%%%%%%%%%%%%%%%%%%%%%%%%%%%%%%%%%%%%%%%%

\noindent
 {\bf Proof.}
As always (see the convention in Section 2), convergence and all relations between random
variables are meant with probability one w.r.t. the measure
$P^\theta$ unless specified otherwise.
 To simplify notation we drop the argument or  the index
 $\theta$ in some of the expressions below.
 Rewrite \eqref{rec2}
in the form
$$
\Dl_t=\Dl_{t-1}+ \Gm_t^{-1}(\theta+\Dl_{t-1})
\p_t(\theta+\Dl_{t-1}).
$$
 By the Taylor expansion,
                           %(liap1)
\begin{eqnarray}
V_t(\Dl_t) =V_t(\Dl_{t-1})+ \dot
V_t(\Dl_{t-1})\Gm_t^{-1}(\theta+\Dl_{t-1})
\p_t(\theta+\Dl_{t-1})  \nonumber\\
 +\frac 12 \left[\Gm_t^{-1}(\theta+\Dl_{t-1})
\p_t(\theta+\Dl_{t-1})\right]^T {\mbox{\"{V}}}_t(\tilde
\Dl_{t}) \Gm_t^{-1}(\theta+\Dl_{t-1}) \p_t(\theta+\Dl_{t-1}),
\nonumber
\end{eqnarray}
Since   $\dot V_t(u)=2u^TC_t$ and ${\mbox{\"{V}}}_t(u)=2C_t$
we obtain
\begin{eqnarray}
V_t(\Dl_t)=V_t(\Dl_{t-1}) + 2\left(C_t \Dl_{t-1},
\Gm_t^{-1}(\theta+\Dl_{t-1})
\p_t(\theta+\Dl_{t-1})\right) \nonumber \\
+ \left[\Gm_t^{-1}(\theta+\Dl_{t-1})
\p_t(\theta+\Dl_{t-1})\right]^T C_t \Gm_t^{-1}(\theta+\Dl_{t-1})
\p_t(\theta+\Dl_{t-1}). \nonumber
\end{eqnarray}
Since
$$
 V_t(\Dl_{t-1})=V_{t-1}(\Dl_{t-1})+\tr V_t(\Dl_{t-1}),
$$
we have
$$
E_\theta \left\{V_t(\Dl_t)\mid {\cf}_{t-1}\right\} =
V_{t-1}(\Dl_{t-1})+ {\cal K}_t.
$$
Then, using the obvious decomposition $ {\cal K}_t= {[{\cal
K}_t]}^+ - {[{\cal K}_t]}^-, $ the previous inequality can be
rewritten as
$$
E_\theta\left\{V_t(\Dl_t)  \mid {\cal{F}}_{t-1}\right\} =
V_{t-1}(\Dl_{t-1})(1+B_t)+B_t- [{\cal K}_t]^-,
$$
where $ B_t=\left(1+V_{t-1}(\Dl_{t-1}) \right)^{-1}[{\cal K}_t]^+.
$ Since, by \eqref{Kt}, $\sum_{t=1}^\infty  B_t < \infty, $
 the assertion of the lemma  follows immediately
 on application of  Lemma A1 in Appendix A
(with $X_n=V_n(\Dl_n)$, $\beta_{n-1}=\xi_{n-1}=B_n$ and
$\zeta_{n-1}={[{\cal K}_n]}^-$).
  $\diamondsuit$
\vskip+0.2cm

%Everywhere below $x^{2\dl}={(x^2)}^\dl$

%\vskip+0.2cm

%%%%%%%%%%%%%%%%%%%%%%%%%%%%%%%%%%%%%%%%%%%%%%%%%%%%%%%%%%%%%%%%%%%%%
%                   C O R O L L A R Y  3.1
%%%%%%%%%%%%%%%%%%%%%%%%%%%%%%%%%%%%%%%%%%%%%%%%%%%%%%%%%%%%%%%%%%%%%
\begin{cor}
Let  $\{a_t(\theta)\}$ be a predictable  non-decreasing scalar process such that
$a_t(\theta)\to \infty$ as $t \to \infty.$ Denote
$\tr a_t(\theta)=a_t(\theta)-a_{t-1}(\theta)$ and suppose that
\begin{description}
\item[(R1)]

$$
 \lim_{t\to \infty}    \frac{\tr a_t(\theta)}{a_{t-1}(\theta)}=0, \qquad
 P^\theta\mbox{-a.s.};
 $$
\item[(R2)] there exist a symmetric and non-negative  definite matrix
$C_\theta$ and a predictable non-negative scalar process
${\cal P}_t$ such that
%                       (R21)
\begin{eqnarray}\label{R21}
2\left( C_\theta \Dl_{t-1}, \Gm_t^{-1}(\theta+\Dl_{t-1})
b_t(\theta,\Dl_{t-1})\right) + {\cal P}_t \leq  \nonumber\\-\lambda_t(\theta)\left(C_\theta
\Dl_{t-1},\Dl_{t-1}\right),
\end{eqnarray}
eventually, where $\{\lambda_t(\theta)\}$ is a predictable scalar process,
satisfying
%                       (R22)
\begin{eqnarray}\label{R22}
\sum_{s=1}^{\infty}{\left[
\frac{\tr a_t(\theta)}{a_t(\theta)}-
\lambda_t(\theta)
\right] }^+ < \infty, \qquad
 P^\theta\mbox{-a.s.};
\end{eqnarray}
 \item[(R3)] for each $0<\ve<1$ and the process ${\cal P}_t$ defined in (R2),
 $$
\sum_{s=1}^{\infty}a_t^{\ve}(\theta)
\left[E_\theta \left\{\|\Gm_t^{-1}(\theta+\Dl_{t-1})
\p_t(\theta+\Dl_{t-1})\|^2 \mid {\cf}_{t-1}\right\}-{\cal P}_t\right]^+
< \infty, \qquad
 P^\theta\mbox{-a.s.}
$$
 \end{description}
Then ~
$
a_t(\theta)^{2\dl}(\hat \theta_t-\theta)^T C_\theta(\hat \theta_t-\theta) \to 0
$
~
 ($P^\theta$-a.s.) for any  $\dl \in ]0,1/2[.$
\end{cor}

%%%%%%%%%%%%%%%%%%%%%%%%%%%%%%%%%%%%%%%%%%%%%%%%%%%%%%%%%%%%%%%%%%%%%
%             P R O O F   OF      C O R      3.1
%%%%%%%%%%%%%%%%%%%%%%%%%%%%%%%%%%%%%%%%%%%%%%%%%%%%%%%%%%%%%%%%%%%%%
\noindent
 {\bf Proof.}  As always (see the convention in Section 2), convergence and all relations between random
variables are meant with probability one w.r.t. the measure
$P^\theta$ unless specified otherwise. Let us check the conditions of Lemma 3.1 for
  $C_t(\theta)=C_\theta(a_t(\theta))^{2\dl},$
     $\dl \in ]0,1/2[$.  To simplify notation we drop the fixed  argument or  the index
 $\theta$ in some of the expressions below.
Denote
$$
r_t={(\tr a_t^{2\dl}- a_t^{2\dl}\lambda_t)}/{a_{t-1}^{2\dl}}
$$
and
$$
\tilde {{\cal P}_t}=a_t^{2\dl}\left({\cal E}_t-{\cal P}_t\right)
$$
where
$$
{\cal E}_t=E_\theta \left\{
\left[\Gm_t^{-1}(\theta+\Dl_{t-1})
\p_t(\theta+\Dl_{t-1})\right]^T C
\left[\Gm_t^{-1}(\theta+\Dl_{t-1}) \p_t(\theta+\Dl_{t-1})\right]
\mid {\cf}_{t-1}\right\}.
$$
By (R2), for  ${\cal K}_t$ defined in
\eqref{Kt2} we have
\begin{eqnarray}
{\cal K}_t
&&=\tr a_t^{2\dl}\left(C
\Dl_{t-1},\Dl_{t-1}\right)+
2a_t^{2\dl}\left( C\Dl_{t-1}, \Gm_t^{-1}(\theta+\Dl_{t-1})
b_t(\theta,\Dl_{t-1})\right)+
\nonumber\\   && ~~~~~~~~~~~~~~~~~~~~~~~~~~~~~~~~~~~~~~~~~~~~~~~~~~~~~~~~~~~~~~~~~~~~~~~~~
a_t^{2\dl}{\cal P}_t +\tilde {\cal P}_t
\nonumber\\
&&\le\left(\tr a_t^{2\dl}- a_t^{2\dl}\lambda_t\right)\left(C
\Dl_{t-1},\Dl_{t-1}\right)+\tilde{{\cal P}_t}
\nonumber\\
&& \le  r_t\left(a_{t-1}^{2\dl}C
\Dl_{t-1},\Dl_{t-1}\right)+\tilde{{\cal P}_t}. \nonumber
\end{eqnarray}
Since $C$ is non-negative definite,
$$
\left(1+ V_{t-1}(\Dl_{t-1})\right)^{-1}\left[
{\cal K}_t\right]^+ =\left(1+ \left(a_{t-1}^{2\dl}C
\Dl_{t-1},\Dl_{t-1}\right)\right)^{-1}\left[
{\cal K}_t\right]^+
\le [r_t]^++[\tilde{{\cal P}_t}]^+.
$$
By (R3),  $\sum_{t=1}^{\infty}[\tilde {{\cal P}_t}]^+ <\infty$  which implies that
 \eqref{Kt} is equivalent to
$
\sum_{t=1}^{\infty}\left[
r_t\right]^+< \infty.
$
Since $\tr a_t^{2\dl}=a_t^{2\dl}-a_{t-1}^{2\dl}$, we can rewrite
$r_t$ as
$$
r_t=
\left( a_ta_{t-1}^{-1}\right)^{2\dl}\left(1-
\lambda_t\right)-1.
$$
Also, since $(1+x)^{2\dl}=1+2\dl x+O(x^2),$ we have
%\begin{equation}
$$
(a_ta_{t-1}^{-1})^{2\dl}=\left(1+\frac {\tr a_t}
{a_{t-1}}\right)^{2\dl}=
1+2\dl\frac {\tr a_t}{a_{t-1}} +\dl_t^{(1)},
$$
%\end{equation}
where, by (R1),
$
\dl_t^{(1)}=O\left( {\tr a_t}/{a_{t-1}}\right)^2 \to 0
$
as $t\to \infty$.
Denote
$$
\eta_t={\tr a_t}/{a_t}-\lambda_t.
$$
Then simple calculations show that
\begin{eqnarray}
r_t
&&\le
\left( a_ta_{t-1}^{-1}\right)^{2\dl}
\left(1+
 {\eta_t}^+-\frac{\tr a_t}{a_t}\right)-1 \nonumber \\
&& =-(1-2\dl)\frac {\tr a_t}{a_{t-1}}
+\dl_t^{(1)}+\eta_t^++2\dl\eta_t^+
\frac {\tr a_t}{a_{t-1}} +
\eta_t^+\dl_t^{(1)}+\nonumber \\&& ~~~~~~~~~~~~~~~~~~~~~~~~~~~~~~~~~
(1-2\dl)\frac{\tr a_t}{a_{t}}
\frac {\tr a_t}{a_{t-1}} -
\frac{\tr a_t}{a_{t}}\dl_t^{(1)} \nonumber \\
&&=\frac {\tr a_t}{a_{t-1}}\left( -(1-2\dl)+\dl_t^{(2)}\right)
+\dl_t^{(3)} \nonumber
\end{eqnarray}
where
$$
\dl_t^{(2)}=\left(\frac {\tr a_t}{a_{t-1}}\right)^{-1}\dl_t^{(1)}
(1-\frac{\tr a_t}{a_{t}})
+(1-2\dl)\frac{\tr a_t}{a_{t}},
$$
$$
\dl_t^{(3)}=\eta_t^++2\dl\eta_t^+
\frac {\tr a_t}{a_{t-1}} +
\eta_t^+\dl_t^{(1)}.
$$
From (R1) and (R2),
$
 \dl_t^{(2)} \to 0 \; \mbox{and}\sum_{t=1}^{\infty}|\dl_t^{(3)}| <
 \infty.
$
Then, since $1-2d>0,$ we obtain that
$
[r_t]^+ \le |\dl_t^{(3)}|.
 $
 It therefore follows that the conditions of Lemma 3.1 are
 satisfied implying that
$ a_t^{2\dl}\|\hat \theta_t-\theta)\|^2 $
 converges to a finite
limit. Finally, since this holds for an arbitrary  $\dl \in
]0,1/2[$
 and  $a_t \to \infty$, the result follows.
$\diamondsuit$

%%%%%%%%%%%%%%%%%%%%%%%%%%%%%%%%%%%%%%%%%%%%%%%%%%%%

%\vskip+0.2cm
%%%%%%%%%%%%%%%%%%%%%%%%%%%%%%%%%%%%%%%%%%%%%%%%%%%%

%%%%%%%%%%%%%%%%%%%%%%%%%%%%%%%%%%%%%%%%%%%%%%%%%%%%%%%%%%%%%%%%%%%%%%%%%%
                           % R E M A R K 3.1
%%%%%%%%%%%%%%%%%%%%%%%%%%%%%%%%%%%%%%%%%%%%%%%%%%%%%%%%%%%%%%%%%%%%%%%%%%
\begin{rem}
{\rm
Note the that the first term in the left hand side of \eqref{R21} is usually negative and assuming
that  ${\cal P}_t=0$ the positive parts in \eqref{R22} are usually zero (or quite small) in many examples.
On the other hand, the choice  ${\cal P}_t=0$
means that (R3) becomes more restrictive  imposing stronger probabilistic
restrictions on the model.
The choice ${\cal P}_t=0$  is natural
in the iid case  since all the required probabilistic conditions are in this case
automatically satisfied.  (see also Remark 3.2).
Now, if the first term in the left hand side of \eqref{R21} is negative  with a
``high enough'' absolute value, then it may be possible to
introduce  a non-zero ${\cal P}_t$  without jeopardising  \eqref{R22}.
One possibility might be  ${\cal P}_t=\|\Gm_t^{-1}(\theta+\Dl_{t-1})b_t(\theta,\Dl_{t-1})\|^2.$
Also, in this case,   since
$b_t(\theta,u)=E_\theta\{
\p_t(\theta+u) \mid {\cf}_{t-1}\}$ and $\Gm_t^{-1}(\theta+u)$ are  predictable processes,
the condition in (R3) can be rewritten as
$$
\sum_{s=1}^{\infty}a_t^{\ve}(\theta)
E_\theta \left\{\|\Gm_t^{-1}(\theta+\Dl_{t-1})\left\{
\p_t(\theta+\Dl_{t-1})-b_t(\theta,\Dl_{t-1})\right\}\|^2 \mid {\cf}_{t-1}\right\}
< \infty.
$$
}
 \end{rem}
%\vskip+0.2cm

%\vskip+0.2cm
%%%%%%%%%%%%%%%%%%%%%%%%%%%%%%%%%%%%%%%%%%%%%%%%%%%%

%%%%%%%%%%%%%%%%%%%%%%%%%%%%%%%%%%%%%%%%%%%%%%%%%%%%%%%%%%%%%%%%%%%%%%%%%%
                           % R E M A R K 3.2
%%%%%%%%%%%%%%%%%%%%%%%%%%%%%%%%%%%%%%%%%%%%%%%%%%%%%%%%%%%%%%%%%%%%%%%%%%
\begin{rem}
{\rm
Consider
the i.i.d. case with
$$
f_ t(\theta, z \mid x_1^{t-1})= f(\theta, z), \;\; \;\;
\p_t(\theta)=\p(\theta, z )|_{z=X_t},
$$
where  $\int \p(\theta,z)f(\theta,z)\mu (dz)=0 $ and $
 \Gm_t(\theta)=t\gm(\theta)
$ for some  invertible  non-random matrix  $\gm(\theta)$. Then
$$
   b_t(\theta, u)=b(\theta, u)=\int
\p(\theta+u,z)f(\theta,z)\mu(\,dz),
$$
implying that $b_t(\theta, 0)=0$. Denote $ \Dl_t=\hat
\theta_t-\theta $ and rewrite \eqref{rec2} in the form
                         % robmonro
\begin{equation}\label{robmonro}
\Dl_t = \Dl_{t-1}  +  \frac1 t\left(\gm^{-1}(\theta+\Dl_{t-1})
b(\theta,\Dl_{t-1})+\ve_t^\theta\right),
\end{equation}
where
$$
\ve_t^\theta=\gm^{-1}(\theta+\Dl_{t-1}) \left\{
\p(\theta+\Dl_{t-1},X_t)-b(\theta,\Dl_{t-1})\right\}.
$$
Equation \eqref{robmonro} defines  a Robbins-Monro stochastic
approximation procedure that converges to  the solution of the
equation
$$
R^{\theta}(u):= \gm^{-1}(\theta+u)b(\theta,u)=0,
$$ when the values of the function
$R^{\theta}(u)$ can only be observed with zero expectation errors
$ \ve_t^\theta$. Note that in  general,  recursion \eqref{rec2}
cannot be considered in the framework of  classical stochastic
approximation theory (see  Lazrieva,  Sharia, and    Toronjadze
    (1997, 2003) for  the generalized   Robbins-Monro
     stochastic approximations procedures).
      For the i.i.d. case, conditions of Corollary 3.1
 can be written as   {\bf (B1)} and {\bf (B2)} in Corollary 4.1 (see also
 Remark 4.1),
 which  are  standard  assumptions for  stochastic approximation
procedures of  type \eqref{robmonro} (see, e.g., Robbins and Monro
(1951), Gladyshev (1965),
 Khas'minskii and Nevelson (1972), Ljung and Soderstrom (1987),   Ljung, Pflug and Walk (1992)).
 }
 \end{rem}
%\vskip+0.2cm

                 % 4  SPECIAL MODELS   AND EXAMPLES

\section{SPECIAL MODELS
                   AND EXAMPLES}
\setcounter{equation}{0}

                      %  E X A M P L E 1
\noindent
{\large \bf 1.} {\bf  The i.i.d. scheme.} Consider  the classical scheme of
  i.i.d. observations $X_1,X_2,\ldots ,$ with a common probability
density/mass  function
$f(\theta,x), \;\; \theta \in {\mathbb{R}}^m.$
Suppose that $\p(\theta, z)$ is an estimating  function with
$$
\int
\p(\theta,z)f(\theta,z)\mu (dz)=0.
$$
Let us define
the recursive estimator $\hat \theta_t$  by
                             % mleg
\begin{equation}\label{mleg}
\hat \theta_t = \hat \theta_{t-1}  + \frac 1 t\gm^{-1}
(\hat \theta_{t-1})
  \p(\hat \theta_{t-1},X_t),\qquad t\ge 1,
\end{equation}
where $\gm(\theta)$ is a non-random matrix such that
$\gm^{-1}(\theta)$ exists  for any $\theta\in {\mathbb{R}}^m$ and
$\hat\theta_0\in {\mathbb{R}}^m$ is any initial value.

%%%%%%%%%%%%%%%%%%%%%%%%%%%%%%%%%%%%%%%%%%%%%%%%%%%%%%%%%%%%%%%%%%%%%%%%%
%             C O R O L L A R Y    4.1
%%%%%%%%%%%%%%%%%%%%%%%%%%%%%%%%%%%%%%%%%%%%%%%%%%%%%%%%%%%%%%%%%%%%%%%%%%
\begin{cor}
Suppose that $\hat \theta\to \theta $ ($P^\theta$-a.s.)  and
\begin{description}
\item[(B1)]
 there exists a symmetric and non-negative  definite matrix $C_\theta$ such that
$$
\left( C_\theta u,\gm^{-1}(\theta+u)E^\theta\p(\theta+u,X_1)
\right) \leq -\frac 12 \left(C_\theta u,u\right),
$$
for small $u$'s;
 \item[(B2)] ~~~~~~~~~~~~~
 $
E_\theta \|\gm^{-1}(\theta+u) \p(\theta+u)\|^2=O(1) $ as $u\to 0$.
\end{description}
Then ~ $ t^\dl(\hat \theta_t-\theta)^T C_\theta (\hat \theta_t-\theta)\to 0$
~ ($P^\theta$-a.s.)  for any  ~ $\dl \in
]0,1/2[.$
\end{cor}
%%%%%%%%%%%%%%%%%%%%%%%%%%%%%%%%%%%%%%%%%%%%%%%%%%%%%%%%%%%%%%%%%%%
                % P R O O F  OF C O R O L L A R Y  4.1
%%%%%%%%%%%%%%%%%%%%%%%%%%%%%%%%%%%%%%%%%%%%%%%%%%%%%%%%%%%%%%%%
 {\bf Proof.} The result follows immediately  if we take $a_t(\theta)=t,$ ${\cal P}_t=0$
 and $\lambda_t(\theta)=1/t$ in Corollary  3.1.
 $\diamondsuit$

%%%%%%%%%%%%%%%%%%%%%%%%%%%%%%%%%%%%%%%%%%%%%%%%%%%%%%%%%%%%%%%%%%%%%%%%%%
                  % R E M A R K 4.1
%%%%%%%%%%%%%%%%%%%%%%%%%%%%%%%%%%%%%%%%%%%%%%%%%%%%%%%%%%%%%%%%%%%%%%%%%%
\begin{rem}
{\rm
 As it was mentioned in
Remark 3.2, for the i.i.d. case the recursive procedures can be
studied in the framework of stochastic approximation theory.
 For stochastic approximation procedures of this type,
conditions which guarantee a good rate of  convergence are
expressed in  terms of stability of matrices. Recall that a matrix
$A$ is called stable if the real parts of its eigenvalues are
negative. A standard  requirement in  stochastic approximation
theory is the existence of the representation (see Remark 3.1 for
the
 notation)
\begin{equation}\label{stab}
R^{\theta}(u)=B^\theta u+o(\|u\|) \;\; \mbox{as} \;\; u \to 0,
\end{equation}
where the matrix $ S^\theta=B^\theta+\frac 12 {{\bf 1}} $ is
stable. It is easy to see that this assumption implies (B1).
 Indeed,
it follows from the stability of $S^\theta$ that the maximum of
the real parts of the eigenvalues of $B^\theta$ is less than ~
$-1/2$.
%for some small $\ve>0$.
This implies (see, e.g., Khas'minskii  and  Nevelson (1972), Ch.6,
\S 3, Corollary 3.1),
 that there
exists   a symmetric and positive definite  matrix $C_\theta$ such
that
$$
\left( C_\theta u, B_\theta u \right) < -\frac 12 \left(C_\theta
u,u\right),
$$
which, together with \eqref{stab}, implies (B1).
}
 \end{rem}

As a particular example, consider
$$
f(\theta, x)=\frac{1}{\pi\left(1+(x-\theta)^2\right)} ,
$$
the probability density function of the Cauchy distribution with mean
$\theta$.
Simple calculations show that
$$
\frac{\dot f}{f}(\theta, x)=\frac{2(x-\theta)}{1+(x-\theta)^2} ~~~\mbox{and}~~~
\frac{\partial^2}{\partial \theta^2}\log f(\theta, x)=\frac{2(x-\theta)^2-2}{(1+(x-\theta)^2)^2}.
$$
Now, using tables of standard integrals, it is easy to check that
$$
i(\theta)=-\int\frac{\partial^2}{\partial \theta^2}\log f(\theta, x)f(\theta, x)\;dx=\frac12.
$$
So, a  ML recursive procedure is
$$
\hat\theta_t=\hat\theta_{t-1}-\frac1t \frac{2(X_t-\hat\theta_{t-1})}{1+(X_t-\hat\theta_{t-1})^2},~~~~~~~~t\ge 1.
$$
Using  tables of standard integrals and  simple algebra,
$$
b(\theta, u)=\frac 2 \pi \int \frac{x-u}{1+(x-u)^2}\frac1{1+x^2}\; dx=-\frac{2u}{4+u^2},
$$
and
$$
\int\left(\frac{\dot f}{f}(\theta+u, x)\right)^2f(\theta, x)\;dx
=\frac 4 \pi \int \left(\frac{x-u}{1+(x-u)^2}\right)^2\frac1{1+x^2}\; dx=\frac{2(4+3u^2)}{(4+u^2)^2}.
$$
Now, it is easy to check that conditions (I) and (II) of Corollary 4.1 in Sharia (2006)
(or in Sharia (1998)) are satisfied, implying that
$
\hat\theta_t \to \theta
$
~~~($P^\theta$-a.s.).
Let us  check the conditions of Corollary 4.1. It follows from
the above calculations that  (B2) holds. Then, for arbitrary  $0<\varepsilon<1/2$
we have
$$
\frac{i^{-1}(\theta)b(\theta,u)}{u}=-\frac{4}{4+u^2}=-1+\frac{u^2}{4+u^2}\le -1 +\varepsilon
$$
for small $u$'s, which yields that (B1) is satisfied with $C_\theta=1.$
Therefore, $ t^\dl(\hat \theta_t-\theta) \to 0$ ~($P^\theta$-a.s.) for any
$0<\delta <1/2.$

%%%%%%%%%%%%%%%%%%%%%%%%%%%%%%%%%%%%%%%%%%%%%%%%%%%%%%%%%%%%%
%      Exponential family of  Markov processes
%%%%%%%%%%%%%%%%%%%%%%%%%%%%%%%%%%%%%%%%%%%%%%%%%%%%%%%%%%%%%
\noindent {\large \bf 2} ~ {\bf  Exponential family of  Markov
processes} Consider a conditional exponential family of  Markov
processes in the sense of Feigin (1981) (see also Barndorf-Nielson
(1988)). This is a time homogeneous Markov chain
with the one-step transition density
$$
f(y; \theta,x)=h(x,y)\exp\left(
\theta^Tm(y,x)-\beta(\theta;x)
\right),
$$
where $m(y,x)$ is a $m$-dimensional vector and $\beta(\theta;x)$
is one dimensional. Then in our notation
$f_t(\theta)=f(X_{t};\theta ,X_{t-1})$ and
$$
l_t(\theta)=\frac{d}{d\theta}\log f_t(\theta)=
m(X_t,X_{t-1})-\dot\beta^T(\theta;X_{t-1}).
$$
 It follows from standard
exponential family theory (see, e.g., Feigin (1981)) that
$l_t(\theta)$
is a martingale-difference and the conditional Fisher information
is
$$
I_t(\theta)=\sum_{s=1}^t\ddot\beta(\theta;X_{s-1}).
$$
So, a maximum likelihood  type recursive procedure  can be defined as
$$
\hat\theta_t=\hat\theta_{t-1}+\left(\sum_{s=1}^t\ddot\beta(\hat\theta_{t-1};X_{s-1})
\right)^{-1}\left(m(X_t,X_{t-1})-\dot\beta^T(\hat\theta_{t-1};X_{t-1})
\right), ~~~ t\ge  1.
$$
Let us find the functions appearing in the conditions of our
theorems for the case $\psi_t=l_t$ and $\Gamma_t=I_t$.
Since
$E_\theta\left\{l_t(\theta)\mid
{\cal{F}}_{t-1}\right\}=0$
we have
$$
E_\theta\left\{m(X_t,X_{t-1}) \mid
{\cal{F}}_{t-1}\right\}=\dot\beta^T(\theta;X_{t-1})
$$
 and also,
\begin{eqnarray}
\ddot\beta(\theta;X_{t-1})=i_t(\theta)=
E_\theta\left\{l_t(\theta)l^T_t(\theta)\mid{\cal{F}}_{t-1}\right\}\nonumber \\
=E_\theta\left\{m(X_t,X_{t-1})m^T(X_t,X_{t-1}) \mid
{\cal{F}}_{t-1}\right\}-\dot\beta^T(\theta;X_{t-1})\dot\beta(\theta;X_{t-1}),
\nonumber
\end{eqnarray}
which implies that
\begin{equation} \label{mm}
E_\theta\left\{m(X_t,X_{t-1})m^T(X_t,X_{t-1}) \mid
{\cal{F}}_{t-1}\right\}=\ddot\beta(\theta;X_{t-1})+
\dot\beta^T(\theta;X_{t-1})\dot\beta(\theta;X_{t-1}).
\end{equation}
Now, it is a simple matter to check that
\begin{eqnarray}\label{bt}
~~~~~~~~b_t(\theta,u)=E_\theta\left\{l_t(\theta+u)\mid
{\cal{F}}_{t-1}\right\}=\dot\beta^T(\theta;X_{t-1})-\dot\beta^T(\theta+u;X_{t-1}).
\end{eqnarray}
Using \eqref{mm} (since ~$\mbox{trace}(vv^T)=v^Tv$ ~and~
$\mbox{trace}(A+B)=$\mbox{trace} A+$\mbox{trace} B$),
\begin{eqnarray}\label{l2}
E_\theta\left\{\|l_t(\theta+u)\|^2\mid
{\cal{F}}_{t-1}\right\}=\mbox{trace}\ddot\beta(\theta;X_{t-1})+
\|\dot\beta^T(\theta;X_{t-1})-\dot\beta^T(\theta+u;X_{t-1})\|^2
\nonumber \\
=\mbox{trace} \ddot\beta(\theta;X_{t-1})+
\|b_t(\theta,u)\|^2.
\end{eqnarray}
Using these expressions one can check conditions of the
relevant theorems  for different choices of  functions $m$  and
$\beta$.

Now suppose that $\theta$ is one dimensional and consider the class of
conditionally additive exponential families, that is,
$$
f(y; \theta,x)=h(x,y)\exp\left(
\theta m(y,x)-\beta(\theta;x)
\right),
$$
with
\begin{eqnarray}\label{Add}
\beta(\theta;x)=\gamma(\theta) h(x)
\end{eqnarray}
where $h(\cdot)\ge 0$  and $\ddot \gamma(\cdot) \ge 0$
(see Feigin (1981)).
Then,
$$
I_t(\theta)=\ddot \gamma(\theta)H_t
~~~ \mbox{where} ~~~
H_t=\sum_{s=1}^t h(X_{s-1}).
$$
Assuming that $\ddot \gamma(\theta)\not= 0,$
the likelihood recursive procedure is
\begin{eqnarray}\label{RecAdd}
\hat\theta_t=\hat\theta_{t-1}+\frac1{\ddot\gamma(\hat\theta_{t-1})H_t}
\left(m(X_t,X_{t-1})-\dot\gamma(\hat\theta_{t-1})h(X_{t-1})
\right).
\end{eqnarray}

The following result gives sufficient conditions for the convergence of
\eqref{RecAdd}.

%%%%%%%%%%%%%%%%%%%%%%%%%%
%      P r o p o s i t i o n  4.1
%%%%%%%%%%%%%%%%%%%%%%%%%%%%%%%%%%%%%%%%%%%%%%%%%%%%%%%%%%%%%%%%%%
\bigskip
\noindent
{\bf Proposition  4.1}
{\it
Suppose that $H_t\to \infty$  ($P^\theta$-a.s.) and either $\dot\gamma$ is a linear function,
or the following conditions are satisfied:
\begin{description}
\item[(M1)]
$$
  \frac{h(X_{t-1})}{H_{t}}\to 0, \qquad
 P^\theta\mbox{-a.s.};
 $$
\item[(M2)]$\;$ for any finite $a$ and $b$,
$$  0 <\inf_{u\in [a,b]} \ddot \gamma(u) \le \sup_{u\in [a,b]} \ddot \gamma(u) <\infty;
$$
\item[(M3)] there exists a constant $B$ such that
$$
\frac{1+\dot\gamma^2(u)}{\ddot\gamma^2(u)} \le B(1+u^2)
$$
for each $ u\in {\mathbb{R}}$.
\end{description}
Then
 $\hat\theta_t$  defined by \eqref{RecAdd} is  strongly consistent (i.e.,
$\hat \theta_t \to \theta \;\; P^\theta$-a.s.) for any  initial
value $\hat \theta_0$ .
}

\medskip
%%%%%%%%%%%%%%%%%%%%%%%%%%%%%%%%%%%%%%%%%%%%%%%%%%%%%%%%%%%%%%%%%%%%%%%%%%
 {\bf Proof.} See Appendix B.

\bigskip

In the next statement we assume that the recursive procedure
converges and study the rate of convergence.

%%%%%%%%%%%%%%%%%%%%%%%%%%%%%%%%%%%%%%%%%%%%%%%%%%%%%%%%%%%%%%%%%%%%%%%%%
%             C O R O L L A R Y    4.2
%%%%%%%%%%%%%%%%%%%%%%%%%%%%%%%%%%%%%%%%%%%%%%%%%%%%%%%%%%%%%%%%%%%%%%%%%%
\begin{cor}
Suppose that $\hat\theta_t$  defined by \eqref{RecAdd} is  strongly consistent (i.e.,
$\hat \theta_t \to \theta \;\; P^\theta$-a.s.). Suppose also that
\begin{description}
\item[(1)] $H_t\to \infty,  \qquad
 P^\theta\mbox{-a.s.}$;
\item[(2)]
$$
  \frac{h(X_{t})}{H_{t}}\to 0, \qquad
 P^\theta\mbox{-a.s.};
 $$
\item[(3)]
$\ddot\gamma(\cdot)$ is a continuous positive  function.
\end{description}
Then ~
$
H_t^\dl(\hat \theta_t-\theta) \to 0$ \;\; ($P^\theta$-a.s.)
~ for any  $\dl \in ]0,1/2[.$
\end{cor}

%%%%%%%%%%%%%%%%%%%%%%%%%%%%%%%%%%%%%%%%%%%%%%%%%%%%%%%%%%%%%%%%%%%%%
%             P R O O F   OF      C O R      4.1
%%%%%%%%%%%%%%%%%%%%%%%%%%%%%%%%%%%%%%%%%%%%%%%%%%%%%%%%%%%%%%%%%%%%%
\noindent
 {\bf Proof.} As always (see the convention in Section 2), convergence and all relations between random
variables are meant with probability one w.r.t. the measure
$P^\theta$ unless specified otherwise.
 By \eqref{bt},
                                % Addbt
\begin{equation}\label{Addbt}
 b_t(\theta,u)=h(X_{t-1})\left(\dot\gamma(\theta)-\dot\gamma(\theta+u)\right).
\end{equation}
Let us check that the conditions of Corollary  3.1 are satisfied with
$\psi_t(\theta)=l_t(\theta)=m(X_t,X_{t-1})-\dot\gamma(\theta)h(X_{t-1}),$  ~ $\Gamma_t(\theta)=I_t(\theta)
=H_t\ddot\gamma(\theta),$ ~
$a_t(\theta)=H_t,$ $C_\theta=1$ and
 ${\cal P}_t=H_t^{-2}\ddot\gamma^{-2}(\theta+\Dl_{t-1})b_t^2(\theta,\Dl_{t-1}).$
Since $\Dl H_t=h(X_{t-1}),$ (R1) is obviously translated into (2).
Since $\dot\gamma(\theta)-\dot\gamma(\theta+u)=-\ddot\gamma(\theta+\tilde u)u$ where
$|\tilde u|\le |u|,$ the left hand side of \eqref{R21} is
$$
-2\frac{h(X_{t-1})}{H_t}\frac{\ddot\gamma(\theta+\tilde\Dl_{t-1})}{\ddot\gamma(\theta+\Dl_{t-1})}\Dl_{t-1}^2
+
\frac{h^2(X_{t-1})}{H_t^2}\left(\frac{\ddot\gamma(\theta+\tilde\Dl_{t-1})}{\ddot\gamma(\theta+\Dl_{t-1})}\right)^2
\Dl_{t-1}^2
$$
Since $\ddot\gamma(\cdot)$ is continuous and $\Dl_{t-1}=\hat \theta_t -\theta \to 0,$
for any small $\tilde\ve >0$ (which may depend on $\theta$),
$1-\tilde\ve <{\ddot\gamma(\theta+\tilde\Dl_{t-1})}/{\ddot\gamma(\theta+\Dl_{t-1})} <1+\tilde\ve$
for large $t$'s.
So,  \eqref{R21} holds with
$$
\lambda_t(\theta)=2(1-\tilde\ve)\frac{h(X_{t-1})}{H_t}-(1+\tilde\ve)^2\frac{h^2(X_{t-1})}{H_t^2}.
$$
To check \eqref{R22}, consider
\begin{eqnarray}\label{Neg}
\frac{h(X_{t-1})}{H_t}-\lambda_t(\theta)=\frac{h(X_{t-1})}{H_t}
\left(
-1+2\tilde\ve+(1+\tilde\ve)^2\frac{h(X_{t-1})}{H_t}
\right)
\end{eqnarray}
Now,  since $\tilde\ve$ is arbitrary,  we can assume
that $-1+2\tilde\ve <0.$
Also,  it follows from (2) that    ${h(X_{t-1})}/{H_t}\to 0$.
Therefore, \eqref{Neg} is negative for large $t$'s, implying that
\eqref{R22} holds true.

To check (R3) note that by  \eqref{l2},
                               % Addp2
\begin{equation}\label{Addp2}
E_\theta \left\{
l_t^2(\theta+u) \mid {\cf}_{t-1}\right\}=\ddot\gamma(\theta)h(X_{t-1})+
b^2_t(\theta,u)
\end{equation}
and so,
$$
H_t^\ve
\left(E_\theta \left\{H_t^{-2}\ddot\gamma^{-2}(\theta+\Dl_{t-1})
l_t^2(\theta+\Dl_{t-1}) \mid {\cf}_{t-1}\right\}-{\cal P}_t
\right)
=\frac{h(X_{t-1})}{H_t^{2-\ve}}\frac{\ddot\gamma(\theta)}{\ddot\gamma^2(\theta+\Dl_{t-1})}.
$$
Now, (R3) follows from (3) and Proposition A2 in Appendix A.
$\diamondsuit$

\bigskip

A particular example of conditional additive  exponential family
is the Gaussian autoregressive model defined by
$$
X_t=\theta X_{t-1}+Z_t, ~~~~~~~~~~~~ t=1,2,\dots,
$$
where $\theta\in{\mathbb{R}}$, $X_0=0$ and $Z_t$'s are independent
random variables with the standard normal distribution.
In this model $m(y,x)=x y$ and
$\beta(\theta,x)=\frac12x^2\theta^2$ so that
we can assume that
$\gamma(\theta)=\theta^2/2$ and $h(x)=x^2.$ Then
$$
l_t(\theta)=X_tX_{t-1}-X^2_{t-1}\theta, ~~~~~I_t=I_t(\theta)
=\sum_{s=1}^tX_{s-1}^2.
$$
Therefore,
                                 % ArLsq
\begin{eqnarray}\label{ArLsq}
&&\hat\theta_t=\hat\theta_{t-1}+\frac1{I_t}
\left(X_tX_{t-1}-X^2_{t-1}\hat\theta_{t-1}\right)\\
&&I_t=I_{t-1}+X^2_{t-1}. \nonumber
\end{eqnarray}
Note that the rate of the conditional Fisher information $I_{t}$ varies for the
different values of $\theta$. Suppose
                                % fish
\begin{equation}\label{fish}
\kappa_t(\theta)=\left\{\begin{array}{lll}
t(1-\theta^2)^{-1}     &\mbox{for $|\theta|<1$} \\
 \frac12 t^2            &\mbox{for $|\theta|=1$} \\
  \theta^{2t}(\theta^2-1)^{-2}                &\mbox{for $|\theta|>1.$}
\end{array}
\right.
\end{equation}
For  $|\theta|<1,$  $I_t/\kappa_t(\theta)\to 1$
in probability as $t\to\infty,$ whereas
$I_t/\kappa_t(\theta)\to W \sim \chi^2(1)$ almost surely
in the  case  $|\theta|>1$ (non-ergodic case).
In the   case $|\theta|=1,$
the ratio $I_t/\kappa_t(\theta)$ converges in distribution, but
not in probability (for details, see White (1958)  and Anderson
(1959)). It is also well known that $I_t\to \infty$  almost surely for any
$\theta \in {\mathbb{R}} $
(see, e.g, Shiryayev (1984), Ch.VII, 5.5).
Also, since $\dot\gamma(\theta)$ is linear  and $H_t=I_t$,  the conditions of Proposition 4.1
are trivially satisfied. Therefore,  for any $\theta \in {\mathbb{R}}, $
the recursive estimator $\hat\theta_t$ is strongly consistent for
any choice of the initial $\hat\theta_0$.

To establish the rate  of convergence we assume that the process is
(strongly) stationary and ergodic. So, $|\theta|<1$ and
and it follows from the ergodic theorem for stationary processes that the limit
                                % ergodic
\begin{equation}\label{ergodic}
\lim_{t\to \infty}\frac 1t I_t
%~~~~~ \mbox{and} ~~~~~
%\lim_{t\to \infty}\frac 1t \sum_{s=1}^t X^4_{s-1}
\end{equation}
exist $P^\theta$-a.s. and is finite
(it can be proved this holds without assumption of strong stationarity.)
Now, taking $H_t=I_t$, we obtain that
$$
\frac {\Delta I_t}{I_{t-1}}=\frac{I_t}{I_{t-1}}-1=\frac t{t-1} d_t-1 \to 0,
$$
since $
d_t=((t-1)/I_{t-1})({I_t}/t)\to 1.
$
This implies  that  (2) of Corollary 4.2  holds. (Note that for the non-ergodic case
$| \theta |>1,$ we do not expect
(2) to hold  since in this case  $\Delta \kappa_t/\kappa_{t-1}=
\theta^2-1 \not\to 0$.)

So, the conditions of Corollary 4.2 are satisfied implying that
$ t^\dl(\hat \theta_t-\theta) \to 0$ for any $0<\delta <1/2.$
\vskip+0.3cm

\bigskip

\begin{center}
{\large \bf APPENDIX A}
\end{center}

%\appendix
%\section{Appendix}
\noindent
{\bf Lemma  A1}
{\it  Let ${\cal F}_0,{\cal F}_1, \dots$ be a non-decreasing sequence of
$\sigma$-algebras and  $X_n, \beta_n, \xi_n, \zeta_n \in {\cal F}_n,
\;\; n\ge 0,$ are nonnegative r.v.'s such that
$$
E(X_n|{\cal F}_{n-1}) \le X_{n-1}(1+\beta_{n-1})+\xi_{n-1}-
\zeta_{n-1}, \;\;\; n\ge 1
$$
eventually. Then}
$$
\{\sum_{i=1}^{\infty} \xi_{i-1}<\infty\}\cap
\{\sum_{i=1}^{\infty} \beta_{i-1} <\infty\}\subseteq
\{X\rightarrow\}\cap\{
\sum_{i=1}^{\infty} \zeta_{i-1} <\infty\}
\quad (P\mbox{-}a.s.),
$$
where $\{X\rightarrow\}$ denotes the set where $\lim_{n\to \infty} X_n$ exists and is
finite.
\vskip+0.2cm
\noindent
{\bf Remark} Proof can be found in
 Robbins and  Siegmund (1971). Note also that this lemma is a  special case of
the theorem on the convergence sets nonnegative semimartingales
(see, e.g., Lazrieva,  Sharia, and  Toronjadze (1997)).

\bigskip
\noindent
{\bf Proposition  A2}
{\it
If  $d_n$ is a  nondecreasing sequence  of positive numbers such that $d_n\to +\infty$,
then
$$
\sum_{n=1}^\infty \tr d_n/d_n=+\infty
$$
 and
$$
\sum_{n=1}^\infty \tr d_n/d_n^{1 + \varepsilon} <+\infty
$$
for any $\varepsilon > 0$.}
\vskip+0.2cm
\noindent
{\bf Proof} The first claim is easily obtained by contradiction from the Kronecker lemma
(see, e.g., Lemma 2, $\S$3, Ch. IV in Shiryayev (1984)). The second one is proved by the
following argument
\begin{eqnarray*}
0 \le \sum_{n=1}^N \frac{\tr d_n}{d_n^{1 + \varepsilon}} \le \sum_{n=1}^N
\int_0^1 \frac{\tr d_n}{(d_{n - 1} + t \tr d_n)^{1 + \varepsilon}}\, dt =
\sum_{n=1}^N \frac{1}{\varepsilon}\left(\frac{1}{d_{n -
1}^\varepsilon} - \frac{1}{d_n^\varepsilon}\right) \\
= \frac{1}{\varepsilon}\left(\frac{1}{d_0^\varepsilon} - \frac{1}{d_N^\varepsilon}\right)
\to \frac{1}{\varepsilon d_0^\varepsilon} <+\infty.
\end{eqnarray*}
$\diamondsuit$

\bigskip

\begin{center}
{\large \bf APPENDIX B}
\end{center}

\noindent
{\bf Theorem B1} (Sharia (2007), Theorem 3.2)
{\em
Suppose that for $\theta \in \mathbb{R}^m$ there exists a real
valued nonnegative function $ V_\theta(u)  :  \mathbb{R}^m
\longrightarrow \mathbb{R} $ having  continuous and bounded
partial second derivatives and
\begin{description}
\item[(G1)] $ V_\theta(0)=0,$ and for each  $\ve\in (0, 1),$
$$
\inf_{ \|u\| \geq \ve} V_\theta(u)> 0 ;
$$
\item[(G2)] there exists a set $A\in \mathcal{F}$  with $P^\theta (A)>0$ such that
\; for each $\ve\in (0, 1),$
$$
\sum_{t=1}^\infty \inf_{\ve \le V_\theta(u) \le {1/\ve}}
\left[{\cal N}_t(u)\right]^-=\infty
$$
\end{description}
on $A$, where
\begin{eqnarray}
{\cal N}_t(u)&= &\dot V_\theta(u) \Gm_t^{-1}(\theta+u)E_\theta \left\{\psi_t(\theta+u)\mid
{\cf}_{t-1}\right\}\nonumber\\
&&
+\frac12\sup_{v} \|{\mbox{\"{V}}_\theta}(v)\| E_\theta
\left\{\|\Gm_t^{-1}(\theta+u) \p_t(\theta+u)\|^2 \mid
{\cf}_{t-1}\right\}, \nonumber
\end{eqnarray}
\;\;\;\;\;\;\;
%  u\in {{\mathbb{R}}}^m;
\begin{description}
\item[(G3)]  $\;$ for $\Dl_t=\hat\theta_t-\theta,$
$$
\sum_{t=1}^\infty (1+V_\theta(\Dl_{t-1}))^{-1} \left[{\cal
N}_t(\Dl_{t-1})\right]^+ < \infty, \qquad
 P^\theta\mbox{-a.s.}.
$$
\end{description}

Then $\hat \theta_t \to \theta \;\; (P^\theta$-a.s.)
 for any  initial value
$\hat \theta_0$, where $\hat \theta_t$ is defined by \ref{rec2}.
}
\medskip

\noindent
{\bf Proof of Proposition 4.2} As always (see the convention in Section 2), convergence and all relations between random
variables are meant with probability one w.r.t. the measure
$P^\theta$ unless specified otherwise.
Let us check that the conditions of Theorem B1 above
are  satisfied with
$\psi_t(\theta)=l_t(\theta)=m(X_t,X_{t-1})-\dot\gamma(\theta)h(X_{t-1}),$  ~ $\Gamma_t(\theta)=I_t(\theta)
=H_t\ddot\gamma(\theta),$ ~ and $V_t=u^2$. Using \eqref{Addbt} and \eqref{Addp2},
we have
                         % AddNt
\begin{eqnarray}\label{AddNt}
{\cal N}_t(u)= 2u\frac 1{H_t\ddot\gamma(\theta+u)}
b_t(\theta,u)+ \frac 1{H_t^2\ddot\gamma^2(\theta+u)}E_\theta
\left\{l_t^2(\theta+u) \mid{\cf}_{t-1}\right\} \nonumber
\\
=\frac {h(X_{t-1})}{H_t}
\frac{\dot\gamma(\theta)-\dot\gamma(\theta+u)}{\ddot\gamma(\theta+u)}u
\left(2+\frac {h(X_{t-1})}{H_t}\frac{\dot\gamma(\theta)-\dot\gamma(\theta+u)}{u\ddot\gamma(\theta+u)}\right)
+\frac {h(X_{t-1})}{H_t^2}\frac{\ddot\gamma(\theta)}{\ddot\gamma^2(\theta+u)}
\nonumber \\
=:{\cal N}_{1t}(u)+{\cal N}_{2t}(u),
\end{eqnarray}
with the convention that $0/0=0.$
Let us show that for large $t$'s,
\begin{equation}\label{Gr1}
2+\frac {h(X_{t-1})}{H_t}\frac{\dot\gamma(\theta)-
\dot\gamma(\theta+u)}{u\ddot\gamma(\theta+u)} \ge 1.
\end{equation}
If $\dot \gamma$ is linear, the above inequality trivially holds since
${h(X_{t-1})}/{H_t}=\Dl H_t/{H_t}\le 1.$ For a non-linear case we have (assuming that $u\not= 0$),
\begin{equation}\label{Pos}
|{(\dot\gamma(\theta)-\dot\gamma(\theta+u))}/{u\ddot\gamma(\theta+u)}|=
\ddot\gamma(\theta+\tilde u)/\ddot\gamma(\theta+u)
\end{equation}
 where $|\tilde u|\le |u|$.
Suppose now that $|u|\le  M$ where $0 <M < \infty.$ Then it
follows from  (M2) that the left hand side of \eqref{Pos} is
bounded by some positive constant. Also, using the obvious inequality
$(a-b)^2\le 2a^2+2b^2$ and   (M3), we obtain that
${(\dot\gamma(\theta)-\dot\gamma(\theta+u))}^2/{\ddot\gamma^2(\theta+u)}\le
\tilde B (1+u^2)$ for any $u$ (where $\tilde B$ may depend on
$\theta$). So, the left hand side of \eqref{Pos} is less than or
equal to $\sqrt{\tilde B(1+u^2)/u^2}=\sqrt{\tilde B(1/u^2+1)}$
 which is bounded by a positive constant if
$|u|\ge  M.$  So, the left hand side of \eqref{Pos} is bounded by a constant
(which may depend on $\theta$) for any $u$. So,  because of (M1) it follows that
 \eqref{Gr1} holds for large $t$'s.
This  implies that ${\cal N}_{1t}(u) \le 0$ for large $t$'s (recall that $\ddot\gamma(\cdot)$ is positive).
So, using (M3) we obtain that for large $t's,$
$$
\frac1{(1+u^2)} \left[{\cal
N}_t(u)\right]^+ \le
\frac1{(1+u^2)}\left[{\cal N}_{2t}(u)\right]^+ \le \frac {h(X_{t-1})}{H_t^2}B_1,
$$
for some constant $B_1$ which may depend on $\theta.$
Now, since
$
\sum_{t=1}^\infty {h(X_{t-1})}/{H_t^2}  < \infty$ (see Proposition A2 in Appendix A),
condition (G3) of Theorem B1
is satisfied. To check condition (G2), note that
$(\dot\gamma(\theta)-\dot\gamma(\theta+u))u \le 0,$
 use  the obvious inequality
$[x]^- \ge -x$, and \eqref{Gr1} to obtain that for large $t$'s
\begin{eqnarray}
\left[{\cal N}_t(u)\right]^- \ge -{\cal N}_{1t}(u)-{\cal N}_{2t}(u) \ge
-\frac {h(X_{t-1})}{H_t}
\frac{\dot\gamma(\theta)-\dot\gamma(\theta+u)}{\ddot\gamma(\theta+u)}u-{\cal N}_{2t}(u)
\nonumber \\
=\frac {h(X_{t-1})}{H_t}
\frac{\ddot\gamma(\theta+\tilde u)}{\ddot\gamma(\theta+u)}u^2-
\frac {h(X_{t-1})}{H_t^2}\frac{\ddot\gamma(\theta)}{\ddot\gamma^2(\theta+u)}
\nonumber
\end{eqnarray}
where $|\tilde u|\le |u|$.
Then, it follows from (M2) that $\sup_{\ve \le |u| \le {1/\ve}} {\ddot\gamma(\theta)}/{\ddot\gamma^2(\theta+u)} <R$
and $\inf_{\ve \le |u| \le {1/\ve}}{\ddot\gamma(\theta+\tilde u)}u^2/
{\ddot\gamma(\theta+u)} >r>0$
(where the positive constants $R$ and $r$
may depend on $\theta$). Note also that these inequalities trivially hold for the linear case.
Therefore, using once more Proposition A2 in Appendix A we obtain that
$
\sum_{t=1}^\infty \inf_{\ve \le |u| \le {1/\ve}}
\left[{\cal N}_t(u)\right]^-=\infty
$
which completes the proof.
$\diamondsuit$

\bigskip

\bigskip
%%%%%%%%%%%%%%%%%%%%%%%%%%%%%%%%%%%%%%%%%%%%%%%%%%%%%%%%%%%%%%%%%%%
%\section                  {R E F E R E N C E S }
%%%%%%%%%%%%%%%%%%%%%%%%%%%%%%%%%%%%%%%%%%%%%%%%%%%%%%%%%%%%%%%%%%%
\begin{center}
       {\Large{\bf  REFERENCES}}
\end{center}
%\vskip+0.5cm
%\begin{enumerate}
%\begin{itemize}
\begin{description}
\item
 {Anderson,} T.W. (1959). On asymptotic distributions of
 estimates of parameters of stochastic difference equations. {\it
 Ann. Math.  Statist.} {\bf 30}, 676--687.
\item
 { Barndorff-Nielsen,} O.E. (1988). {\it Parametric Statistical
 Models and Likelihood.} Springer Lecture Notes in Statistics 50.
 Heidelberg, Springer.
\item
 {Barndorff-Nielsen,} O.E. and {Sorensen,} M. (1994).
 A review of some aspects of asymptotic likelihood theory for
 stochastic processes. {\it International Statistical Review.}
 {\bf 62,} 1, 133-165.
\item
     {Campbell,} K. (1982). Recursive computation of
     M-estimates for the parameters of a finite autoregressive
     process. {\it Ann. Statist.} {\bf 10}, 442-453.
\item
 { Englund,} J.-E.,   {Holst,} U., and {Ruppert,} D. (1989)
  Recursive estimators for
         stationary, strong mixing processes -- a representation
         theorem and asymptotic distributions
         {\it Stochastic Processes Appl.}  {\bf 31}, 203--222.
\item
{Fabian,} V. (1978).  On asymptotically  efficient
         re\-cur\-sive es\-ti\-ma\-tion,
         {\it Ann. Statist.} {\bf 6}, 854-867.
 \item
{Feigin,} P.D. (1981).  Conditional exponential families and
a representation theorem for asymptotic inference.
{\it Ann. Statist.} {\bf 9}, 597-603.
\item
 {Gladyshev,} E.G. (1965). {On stochastic approximation.}
  {\it Theory Probab. Appl.} {\bf 10}, 297--300.
\item
  {Hampel,} F.R.,   {Ronchetti,} E.M.,   {Rousseeuw,}
  P.J., and {\sc Stahel,} W. (1986).
         {\it Robust Statistics - The Approach Based on Influence
         Functions.} Wiley, New York
   {Huber,} P.J. (1981). {\it Robust Statistics.} Wiley, New York.

\item
   {Jure$\check{\mbox{{c}}}$kov$\acute{\mbox{{a}}}$,} J.  and
    {Sen,} P.K. (1996). {\it Robust Statistical Procedures -
    Asymptotics and Interrelations.} Wiley, New York.
\item
  {Khas'minskii,} R.Z. and {Nevelson,} M.B.     (1972).
  {\it Stochastic Approximation and
          Recursive Estimation.} Nauka, Moscow.
\item
   {Launer,} R.L.  and  {Wilkinson,}  G.N. (1979).
   {\it Robustness in  Statistics.} Academic Press, New York.
\item
    {Lazrieva,}  N., {Sharia,} T. and  {Toronjadze,} T.
    (1997). The Robbins-Monro type
     stochastic differential equations. I. Convergence of solutions.
      {\it Stochastics and Stochastic Reports} {\bf 61}, 67--87.
\item
    {Lazrieva,}  N., {Sharia,} T. and  {Toronjadze,} T.
    (2003). The Robbins-Monro type
     stochastic differential equations. II.  Asymptotic behaviour of solutions.
      {\it Stochastics and Stochastic Reports} {\bf 75}, 153--180.
\item
    {Lazrieva,} N. and  {Toronjadze,} T. (1987).  Ito-Ventzel's formula for
semimartingales, asymptotic properties of MLE and recursive
    estimation.  {\it Lect. Notes in Control and Inform. Sciences, 96, Stochast.
    diff.  systems, H.J, Engelbert, W. Schmidt (Eds.), Springer}
  346--355.
\item
     {Lehmann,} E.L. (1983). {\it Theory of Point Estimation.} Wiley, New
      York.
\item
 {Ljung,} L.  {Pflug,} G. and {Walk,} H.  (1992).
 {\it Stochastic Approximation and
         Optimization of Random Systems.}  Birkh\"auser, Basel.
\item  {Ljung,} L. and {Soderstrom,} T. (1987). {\it Theory and
Practice of Recursive Identification.} MIT Press.
\item
        {Prakasa Rao,} B.L.S. (1999). {\it
        Semimartingales and their Statistical Inference.}
     Chapman $\&$ Hall,   New York.
\item  {Rieder,} H. (1994). {\it Robust Asymptotic Statistics.}
         Springer--Verlag,  New York.
\item
  {Robbins,} H. and {Monro,} S. (1951) A stochastic approximation method,
         {\it Ann. Statist.}  {\bf 22},  400--407.
\item
    {Robbins,} H.  and {Siegmund,}  D. (1971) A convergence theorem for
          nonnegative almost supermartingales and some applications,
         {\it Optimizing  Methods in Statistics}, ed. J.S. Rustagi
          Academic Press, New York.  233--257.
\item
 {Serfling,} R.J. (1980). {\it Approximation Theorems of
 Mathematical  Statistics.} Wiley, New York.
%\item
 %{Sharia,} T.   (1992).
  %On the recursive parameter estimation for general statistical model
%in discrete time.
 %Bulletin of the Georgian Acad. Scienc. 145  {\bf 3}, 465--468.
\item
    {Sharia,} T. (1998). On the recursive parameter estimation for the
          general discrete time statistical model. {\it
          Stochastic Processes Appl.}
          {\bf 73}, {\bf 2}, 151--172.
\item
 {Sharia,} T.   (1997).
 Truncated recursive estimation
 procedures.
          {\it Proc. A. Razmadze Math. Inst.} {\bf  115}, 149--159.
\item
 {Sharia,} T.   (2007).
 Recursive parameter estimation: convergence. {\it Statistical Inference for Stochastic Processes}
 (to appear).
(See also \\ {\em  http://personal.rhul.ac.uk/UkAH/113/GmjA.pdf}).
 \item
  {Shiryayev,} A.N. (1984). {\it Probability.}
           Springer-Verlag, New York.
   \item
 {White,} J.S. (1958). The limiting distribution of the serial
 correlation coefficient in the explosive case. {\it Ann. Math. Stat.}
 {\bf 29}, 1188--1197.
\end{description}
%\end{enumerate}
%\end{itemize}

\end{document}